\documentclass[12pt]{article}
\usepackage{amsfonts}
\usepackage{full page,
amssymb, amscd, amsmath,graphicx, 
}

 \title{{\bf Differential equations and logarithmic intertwining operators for strongly graded vertex algebras}}
 \author{Jinwei Yang}
    \date{}
    \begin{document}
    \bibliographystyle{alpha}
    \maketitle
\newtheorem{thm}{Theorem}[section]
\newtheorem{defn}[thm]{Definition}
\newtheorem{prop}[thm]{Proposition}
\newtheorem{cor}[thm]{Corollary}
\newtheorem{lemma}[thm]{Lemma}
\newtheorem{rema}[thm]{Remark}
\newtheorem{app}[thm]{Application}
\newtheorem{prob}[thm]{Problem}
\newtheorem{conv}[thm]{Convention}
\newtheorem{conj}[thm]{Conjecture}
\newtheorem{cond}[thm]{Condition}
    \newtheorem{exam}[thm]{Example}
\newtheorem{assum}[thm]{Assumption}
     \newtheorem{nota}[thm]{Notation}
\newcommand{\halmos}{\rule{1ex}{1.4ex}}
\newcommand{\pfbox}{\hspace*{\fill}\mbox{$\halmos$}}

\newcommand{\nn}{\nonumber \\}

 \newcommand{\res}{\mbox{\rm Res}}
 \newcommand{\ord}{\mbox{\rm ord}}
\renewcommand{\hom}{\mbox{\rm Hom}}
\newcommand{\edo}{\mbox{\rm End}\ }
 \newcommand{\pf}{{\it Proof.}\hspace{2ex}}
 \newcommand{\epf}{\hspace*{\fill}\mbox{$\halmos$}}
 \newcommand{\epfv}{\hspace*{\fill}\mbox{$\halmos$}\vspace{1em}}
 \newcommand{\epfe}{\hspace{2em}\halmos}
  \newcommand{\nno}{\nonumber}
\newcommand{\nord}{\mbox{\scriptsize ${\circ\atop\circ}$}}
\newcommand{\wt}{\mbox{\rm wt}\ }
\newcommand{\swt}{\mbox{\rm {\scriptsize wt}}\ }
\newcommand{\lwt}{\mbox{\rm wt}^{L}\;}
\newcommand{\rwt}{\mbox{\rm wt}^{R}\;}
\newcommand{\slwt}{\mbox{\rm {\scriptsize wt}}^{L}\,}
\newcommand{\srwt}{\mbox{\rm {\scriptsize wt}}^{R}\,}
\newcommand{\clr}{\mbox{\rm clr}\ }
\newcommand{\tr}{\mbox{\rm Tr}}
\newcommand{\C}{\mathbb{C}}
\newcommand{\Z}{\mathbb{Z}}
\newcommand{\R}{\mathbb{R}}
\newcommand{\Q}{\mathbb{Q}}
\newcommand{\N}{\mathbb{N}}
\newcommand{\CN}{\mathcal{N}}
\newcommand{\F}{\mathcal{F}}
\newcommand{\I}{\mathcal{I}}
\newcommand{\V}{\mathcal{V}}
\newcommand{\one}{\mathbf{1}}
\newcommand{\BY}{\mathbb{Y}}
\newcommand{\ds}{\displaystyle}

        \newcommand{\ba}{\begin{array}}
        \newcommand{\ea}{\end{array}}
        \newcommand{\be}{\begin{equation}}
        \newcommand{\ee}{\end{equation}}
        \newcommand{\bea}{\begin{eqnarray}}
        \newcommand{\eea}{\end{eqnarray}}
         \newcommand{\lbar}{\bigg\vert}
        \newcommand{\p}{\partial}
        \newcommand{\dps}{\displaystyle}
        \newcommand{\bra}{\langle}
        \newcommand{\ket}{\rangle}

        \newcommand{\ob}{{\rm ob}\,}
        \renewcommand{\hom}{{\rm Hom}}

\newcommand{\A}{\mathcal{A}}
\newcommand{\Y}{\mathcal{Y}}

\newcommand{\dlt}[3]{#1 ^{-1}\delta \bigg( \frac{#2 #3 }{#1 }\bigg) }

\newcommand{\dlti}[3]{#1 \delta \bigg( \frac{#2 #3 }{#1 ^{-1}}\bigg) }

 \makeatletter
\newlength{\@pxlwd} \newlength{\@rulewd} \newlength{\@pxlht}
\catcode`.=\active \catcode`B=\active \catcode`:=\active
\catcode`|=\active
\def\sprite#1(#2,#3)[#4,#5]{
   \edef\@sprbox{\expandafter\@cdr\string#1\@nil @box}
   \expandafter\newsavebox\csname\@sprbox\endcsname
   \edef#1{\expandafter\usebox\csname\@sprbox\endcsname}
   \expandafter\setbox\csname\@sprbox\endcsname =\hbox\bgroup
   \vbox\bgroup
  \catcode`.=\active\catcode`B=\active\catcode`:=\active\catcode`|=\active
      \@pxlwd=#4 \divide\@pxlwd by #3 \@rulewd=\@pxlwd
      \@pxlht=#5 \divide\@pxlht by #2
      \def .{\hskip \@pxlwd \ignorespaces}
      \def B{\@ifnextchar B{\advance\@rulewd by \@pxlwd}{\vrule
         height \@pxlht width \@rulewd depth 0 pt \@rulewd=\@pxlwd}}
      \def :{\hbox\bgroup\vrule height \@pxlht width 0pt depth
0pt\ignorespaces}
      \def |{\vrule height \@pxlht width 0pt depth 0pt\egroup
         \prevdepth= -1000 pt}
   }
\def\endsprite{\egroup\egroup}
\catcode`.=12 \catcode`B=11 \catcode`:=12 \catcode`|=12\relax
\makeatother

\def\hboxtr{\FormOfHboxtr} 
\sprite{\FormOfHboxtr}(25,25)[0.5 em, 1.2 ex] 

:BBBBBBBBBBBBBBBBBBBBBBBBB | :BB......................B |
:B.B.....................B | :B..B....................B |
:B...B...................B | :B....B..................B |
:B.....B.................B | :B......B................B |
:B.......B...............B | :B........B..............B |
:B.........B.............B | :B..........B............B |
:B...........B...........B | :B............B..........B |
:B.............B.........B | :B..............B........B |
:B...............B.......B | :B................B......B |
:B.................B.....B | :B..................B....B |
:B...................B...B | :B....................B..B |
:B.....................B.B | :B......................BB |
:BBBBBBBBBBBBBBBBBBBBBBBBB |

\endsprite
\def\shboxtr{\FormOfShboxtr} 
\sprite{\FormOfShboxtr}(25,25)[0.3 em, 0.72 ex] 

:BBBBBBBBBBBBBBBBBBBBBBBBB | :BB......................B |
:B.B.....................B | :B..B....................B |
:B...B...................B | :B....B..................B |
:B.....B.................B | :B......B................B |
:B.......B...............B | :B........B..............B |
:B.........B.............B | :B..........B............B |
:B...........B...........B | :B............B..........B |
:B.............B.........B | :B..............B........B |
:B...............B.......B | :B................B......B |
:B.................B.....B | :B..................B....B |
:B...................B...B | :B....................B..B |
:B.....................B.B | :B......................BB |
:BBBBBBBBBBBBBBBBBBBBBBBBB |

\endsprite

\vspace{2em}



\renewcommand{\theequation}{\thesection.\arabic{equation}}
\renewcommand{\thethm}{\thesection.\arabic{thm}}
\setcounter{equation}{0} \setcounter{thm}{0} 
\date{}
\maketitle

\begin{abstract}
We derive certain systems of differential equations for matrix elements of products and iterates of logarithmic intertwining operators among strongly graded generalized modules for a strongly graded vertex algebra under a certain finiteness condition and a condition related to the horizontal gradings. Using these systems of differential equations, we verify the convergence and extension property needed in the logarithmic tensor category theory for such strongly graded generalized modules developed by Huang, Lepowsky and Zhang.
\end{abstract}

\section{Introduction}
In the present paper, we generalize the arguments in \cite{H} (cf. \cite{HLZ7}) to prove that for a {\it strongly graded} conformal vertex algebra $V$, matrix elements of products and iterates of {\it logarithmic} intertwining operators among triples of strongly graded generalized $V$-modules under suitable assumptions satisfy certain systems of differential equations and that the prescribed singular points are regular. Using these differential equations, we verify the {\it convergence and extension property} needed in the theory of logarithmic tensor categories for strongly graded generalized $V$-modules in \cite{HLZ7}.

The notion of strongly graded conformal vertex algebra and the notion of its strongly graded module were
introduced in \cite{HLZ1} as natural concepts from which the theory of logarithmic tensor categories was developed. A strongly $A$-graded conformal vertex algebra $V$ (respectively, a strongly $\tilde{A}$-graded $V$-module) is a vertex algebra (respectively, a $V$-module), with a weight-grading provided by a conformal vector in $V$ (an $L(0)$-eigenspace decomposition), and with a second, compatible grading by an abelian group $A$ (respectively, an abelian group $\tilde{A}$ containing $A$ as its subgroup), satisfying certain grading restriction conditions. One important source of examples of strongly graded conformal vertex algebras and modules comes {}from the vertex algebras and modules associated with {\it not necessarily positive definite} even lattices. In particular, the tensor products of vertex operator algebras and the vertex algebras associated with even lattices are strongly graded conformal vertex algebras (see \cite{Y1},\cite{Y2}). In \cite{B1}, Borcherds used the vertex algebra associated with the self-dual Lorentzian lattice of rank $2$ and its tensor product with $V^{\natural}$ to construct the ``Monster" Lie algebra.

It was proved in \cite{H} that if every module $W$ for a vertex operator algebra $V = \coprod_{n \in \mathbb{Z}} V_{(n)}$ satisfies the {\em $C_1$-cofiniteness condition}, that is, dim $W/C_1(W) < \infty$, where $C_1(W)$ is the subspace of $W$ spanned by elements of the form $u_{-1}w$ for $u \in V_{+} = \coprod_{n > 0}V_{(n)}$ and $w \in W$, then matrix elements of products and iterates of intertwining operators among triples of $V$-modules satisfy certain systems of differential equations. Moreover, for prescribed singular points, there exist such systems of differential equations such that the prescribed singular points are regular.

To develop the representation theory of vertex operator algebras that are not reductive,
it is necessary to consider generalized modules that are not completely reducible and
the logarithmic intertwining operators among them (see \cite{M1} and \cite{M2}). In \cite{HLZ7}, using the same argument as in \cite{H}, certain systems of differential equations were derived for matrix elements of products and iterates of logarithmic intertwining operators among triples of generalized $V$-modules. In this paper, we prove similar, more general results for matrix elements of products and iterates of logarithmic intertwining operators among triples of strongly graded generalized modules for a strongly graded vertex algebra.

We first generalize the $C_1$-cofiniteness condition for generalized modules for a vertex operator algebra to a {\it $C_1$-cofiniteness condition with respect to $\tilde{A}$} for strongly $\tilde{A}$-graded generalized modules for a strongly graded vertex algebra. That is, every strongly graded generalized $\tilde{A}$-module $W$ satisfies the condition that for $\beta \in \tilde{A}$, dim $W^{(\beta)}/(C_1(W))^{(\beta)} < \infty$, where $W^{(\beta)}$ and $(C_1(W))^{(\beta)}$ are the $\tilde{A}$-homogeneous subspace of $W$ and $C_1(W)$ with $\tilde{A}$-grading $\beta$, respectively. Furthermore, we associate to each $w \in W^{(\beta)}$ a set of partitions of $\beta$ in $\tilde{A}$, which we call an $\tilde{A}$-{\it pattern} of $w$. We say vertex operators preserve the $\tilde{A}$-pattern of an element $w \in W^{(\beta)}$ if for each $A$-homogeneous element $u \in V^{(\alpha)}$, $k \in \Z$ and the $\tilde{A}$-pattern of $w$, the $\tilde{A}$-pattern of $u_kw$ only involve $\alpha$ and the $\tilde{A}$-pattern of $w$. This is a very natural assumption and is very easy to verify for familiar examples of strongly graded modules for strongly graded vertex algebras.

The key step in deriving systems of differential equations in \cite{H} is to construct an $R = \C[z_1^{\pm 1}, z_2^{\pm 1}, (z_1 - z_2)^{-1}]$-module $T$, which is a tensor product of $R$ and a given quadruple of modules for a vertex operator algebra, and a finitely generated quotient module of $T$ by dividing an $R$-submodule $J$ equivalent to the Jacobi identity of intertwining operators. However, for a quadruple of strongly graded generalized modules for a strongly graded vertex algebra, both the $R$-module $T$ constructed in the same way involves tensor products of infinitely many $\tilde{A}$-homogeneous subspaces of the strongly graded generalized modules and thus the quotient module $T/J$ is not finitely generated.

Thanks to the assumption that vertex operators preserve the $\tilde{A}$-patterns of elements in the $C_1$-cofinite strongly graded generalized modules, we can consider a submodule of $T$ consisting of tensor products of four elements whose $\tilde{A}$-patterns come from a given quadruple of elements in the quadruple of strongly graded generalized modules. In particular, this $R$-submodule of $T$ only involves finitely many $\tilde{A}$-homogeneous subspaces and its quotient module by $J$ is finitely generated.

Under the $C_1$-cofiniteness condition with respect to $\tilde{A}$ and the assumption that the $\tilde{A}$-pattern of each element is preserved under vertex operators for the strongly $\tilde{A}$-graded generalized modules, we construct a natural map from the finitely generated $R$-module to the set of matrix elements of products and iterates of logarithmic intertwining operators among triples of strongly graded generalized $V$-modules. The images of certain elements under this map provide systems of differential equations for the matrix elements of products and iterates of logarithmic intertwining operators, as a consequence of the $L(-1)$-derivative property for the logarithmic intertwining operators. Moreover, for any prescribed singular point, we derive certain systems of differential equations such that this prescribed singular point is regular. Using these systems of differential equations, we verify the convergence and extension property needed in the construction of associativity isomorphism for the logarithmic tensor category structure developed in \cite{HLZ1}-\cite{HLZ8}.

The present paper is organized as follows: In Section $2$, we recall the definitions and some basic properties of strongly graded vertex algebras and their strongly graded generalized modules. In Section 3, we introduce the $C_1$-cofiniteness condition with respect to $\tilde{A}$ for strongly $\tilde{A}$-graded generalized modules, the $\tilde{A}$-patterns of elements in the $C_1$-cofinite strongly graded generalized modules and the assumption that the $\tilde{A}$-patterns are preserved by the vertex operators. In Section $4$, we recall the definition of logarithmic intertwining operators among strongly graded generalized modules. The existence of systems of differential equations and the existence of systems with regular prescribed singular points are established in Sections $5$ and $6$, respectively. In Section $7$, we prove the convergence and extension property for products and iterates of logarithmic intertwining operators among strongly graded generalized modules for a strongly graded vertex algebra. In Section $8$, we provide examples of strongly graded vertex algebras and their strongly graded generalized modules.\\

\paragraph{Acknowledgments}
I would like to thank Professor Yi-Zhi Huang and Professor James Lepowsky for helpful discussions and suggestions. I am also grateful for partial support from NSF grant PHY-0901237.

\setcounter{equation}{0}
\section{Strongly graded vertex algebras and their modules}
In this Section, we recall the basic definitions from \cite{HLZ1} (cf. \cite{Y1}, \cite{Y2}).
\begin{defn}
{\rm A {\it conformal vertex algebra} is a ${\mathbb Z}$-graded
vector space
\[
V=\coprod_{n\in{\mathbb Z}} V_{(n)}
\]
equipped with a linear map:
\begin{eqnarray*}
V&\to&({\rm End}\; V)[[x, x^{-1}]] \nn
v&\mapsto& Y(v, x)={\displaystyle \sum_{n\in{\mathbb
Z}}}v_{n}x^{-n-1},
\end{eqnarray*}
and equipped also with two distinguished vectors: {\it vacuum vector} ${\bf 1}\in
V_{(0)}$ and {\it conformal vector} $\omega\in V_{(2)}$, satisfying the following conditions for $u,v \in
V$:
\begin{itemize}
\item the {\it lower truncation condition}:
\[
u_{n}v=0\;\;\mbox{\; for} \;\; n \mbox{ \;sufficiently\; large};
\]
\item the {\it vacuum property}:
\[
Y({\bf 1}, x)=1_V;
\]
\item the {\it creation property}:
\[
Y(v, x){\bf 1} \in V[[x]]\;\;\mbox{ and }\;\lim_{x\rightarrow 0}Y(v,
x){\bf 1}=v;
\]
\item the {\it Jacobi identity} (the
main axiom):
\begin{eqnarray*}
&x_0^{-1}\delta \bigg({\displaystyle\frac{x_1-x_2}{x_0}}\bigg)Y(u,
x_1)Y(v, x_2)-x_0^{-1} \delta
\bigg({\displaystyle\frac{x_2-x_1}{-x_0}}\bigg)Y(v, x_2)Y(u,
x_1)&\nno \\ &=x_2^{-1} \delta
\bigg({\displaystyle\frac{x_1-x_0}{x_2}}\bigg)Y(Y(u, x_0)v,
x_2);
\end{eqnarray*}
\item the {\em Virasoro algebra relations}:
\[
[L(m), L(n)]=(m-n)L(m+n)+{\displaystyle\frac{1}{12}}
(m^3-m)\delta_{n+m,0}c
\]
for $m, n \in {\mathbb Z}$, where
\[
L(n)=\omega _{n+1}\;\; \mbox{ for } \;  n\in{\mathbb Z}, \;\;
\mbox{i.e.},\;\;Y(\omega, x)=\sum_{n\in{\mathbb Z}}L(n)x^{-n-2},
\]
\[
c\in {\mathbb C}\;\;\; (\mbox{central charge of}\; V);
\]
satisfying the {\it  $L(-1)$-derivative property}:
\[
{\displaystyle \frac{d}{dx}}Y(v, x)=Y(L(-1)v, x);
\]
and
\[
L(0)v=nv=(\wt v)v \;\; \mbox{ for }\; n\in {\mathbb Z}\; \mbox{ and
}\; v\in V_{(n)}.
\]
\end{itemize}
}
\end{defn}

This completes the definition of the notion of conformal vertex
algebra.  We will denote such a conformal vertex algebra by
$(V,Y,{\bf 1},\omega)$. \\

\begin{defn}\label{cvamodule}
{\rm  Given a conformal vertex algebra $(V,Y,{\bf 1},\omega)$,  a
{\it module} for $V$ is a ${\mathbb C}$-graded vector space
\begin{equation}\label{grading}
W=\coprod_{n\in{\mathbb C}} W_{(n)}
\end{equation}
equipped with a linear map
\begin{eqnarray*}
V &\rightarrow & (\mbox{End}\ W)[[x,x^{-1}]] \nn
v&\mapsto & Y(v,x) =\sum_{n\in {\mathbb Z}}v_nx^{-n-1}
\end{eqnarray*}
such that the following conditions are
satisfied:
\begin{itemize}
\item
the lower truncation condition: for $v \in V$ and $w \in
W$,
\[
v_nw = 0 \;\;\mbox{ for }\;n\;\mbox{ sufficiently large};
\]
\item
the vacuum property:
\[
Y(\mbox{\bf 1},x) = 1_W;
\]
\item
the Jacobi identity for vertex operators on $W$: for $u, v \in V$,
\begin{eqnarray*}
&{\dps x^{-1}_0\delta \bigg( {x_1-x_2\over x_0}\bigg)
Y(u,x_1)Y(v,x_2) - x^{-1}_0\delta \bigg( {x_2-x_1\over -x_0}\bigg)
Y(v,x_2)Y(u,x_1)}&\nno\\
&{\dps = x^{-1}_2\delta \bigg( {x_1-x_0\over x_2}\bigg)
Y(Y(u,x_0)v,x_2)};
\end{eqnarray*}
\item
the Virasoro algebra relations on $W$ with
scalar $c$ equal to the central charge of $V$:
\[
[L(m), L(n)]=(m-n)L(m+n)+{\displaystyle\frac{1}{12}}
(m^3-m)\delta_{n+m,0}c
\]
for $m,n \in {\mathbb Z}$, where
\[
L(n)=\omega _{n+1}\;\; \mbox{ for }n\in{\mathbb Z}, \;\;{\rm
i.e.},\;\;Y(\omega, x)=\sum_{n\in{\mathbb Z}}L(n)x^{-n-2};
\]
satisfying the $L(-1)$-derivative property
\[
\displaystyle \frac{d}{dx}Y(v, x)=Y(L(-1)v, x);
\]
and
\begin{equation}\label{weight}
(L(0)-n)w=0\;\;\mbox{ for }\;n\in {\mathbb C}\;\mbox{ and }\;w\in
W_{(n)}.
\end{equation}
\end{itemize}
}
\end{defn}

This completes the definition of the notion of module for a
conformal vertex algebra.

\begin{defn}{\rm A generalized module for a conformal vertex algebra is defined in the same way as a module for a conformal vertex algebra except that in the grading (\ref{grading}), each space $W_{(n)}$ is replaced by $W_{[n]}$, where $W_{[n]}$ is the generalized $L(0)$-eigenspace corresponding to the generalized eigenvalue $n \in \C$; that is, (\ref{grading}) and (\ref{weight}) in the definition are replaced by
\[
W=\coprod_{n\in{\mathbb C}} W_{[n]}
\]
and
\[
{\rm for}\ n \in \C \ {\rm and}\ w \in W_{[n]},\ (L(0) - n)^kw = 0, \ {\rm for}\ k \in \N\ {\rm sufficiently\ large},
\]
respectively. For $w \in W_{[n]}$, we still write wt $w = n$ for the generalized weight of $w$.}
\end{defn}

\begin{defn}\label{def:dgv}
{\rm Let $A$ be an abelian group.  A conformal vertex algebra
\[
V=\coprod_{n\in {\mathbb Z}} V_{(n)}
\]
is said to be {\em strongly graded with respect to $A$} (or {\em
strongly $A$-graded}, or just {\em strongly graded} if the abelian
group $A$ is understood) if it is equipped with a second gradation, by $A$,
\[
V=\coprod _{\alpha \in A} V^{(\alpha)},
\]
such that the following conditions are satisfied: the two gradations
are compatible, that is,
\[
V^{(\alpha)}=\coprod_{n\in {\mathbb Z}} V^{(\alpha)}_{(n)}, \;\;
\mbox{where}\;V^{(\alpha)}_{(n)}=V_{(n)}\cap V^{(\alpha)}\;
\mbox{ for any }\;\alpha \in A;
\]
for any $\alpha,\beta\in A$ and $n\in {\mathbb Z}$,
\begin{eqnarray*}
&&V^{(\alpha)}_{(n)}=0\;\;\mbox{ for }\;n\;\mbox{ sufficiently
negative};\\
&&\dim V^{(\alpha)}_{(n)} <\infty;\\
&&{\bf 1}\in V^{(0)}_{(0)};\;\;\;\;\omega\in V^{(0)}_{(2)};\\
&&v_l V^{(\beta)} \subset V^{(\alpha+\beta)}\;\; \mbox{ for any
}\;v\in V^{(\alpha)},\;l\in {\mathbb Z}.
\end{eqnarray*}
}
\end{defn}

This completes the definition of the notion of strongly $A$-graded
conformal vertex algebra.\\

For modules for a strongly graded algebra we would also like to have a second
grading by an abelian group, and it is natural to allow this group
to be larger than the second grading group $A$ for the algebra.
(Note that this already occurs for the {\em first} grading group,
which is ${\mathbb Z}$ for algebras and ${\mathbb C}$ for modules.)

\begin{defn}\label{def:dgw}{\rm
Let $A$ be an abelian group and $V$ a strongly $A$-graded conformal
vertex algebra. Let $\tilde A$ be an abelian group containing $A$ as
a subgroup. A $V$-module (respectively, generalized $V$-module)
\[
W=\coprod_{n\in{\mathbb C}} W_{(n)}\; \;(\mbox{respectively, }\; W^{(\beta)}=\coprod_{n\in {\mathbb C}}W_{[n]})
\]
is said to be {\em strongly graded with respect to $\tilde A$} (or
{\em strongly $\tilde A$-graded}, or just {\em strongly graded}) if
the abelian group $\tilde A$ is understood) if it is equipped with a
second gradation, by $\tilde A$,
\[
W=\coprod _{\beta \in \tilde A} W^{(\beta)},
\]
such that the following conditions are satisfied: the two gradations
are compatible, that is, for any $\beta \in \tilde A$,
\[
W^{(\beta)}=\coprod_{n\in {\mathbb C}} W^{(\beta)}_{(n)},
\;\;\mbox{where }\; W^{(\beta)}_{(n)}=W_{(n)}\cap W^{(\beta)}
\]
\[
(\mbox{respectively, }\; W^{(\beta)}=\coprod_{n\in {\mathbb C}}
W^{(\beta)}_{[n]}, \;\;\mbox{where }\;
W^{(\beta)}_{[n]}=W_{[n]}\cap W^{(\beta)});
\]
for any $\alpha\in A$, $\beta\in \tilde A$ and $n\in {\mathbb C}$,
\begin{eqnarray}
&&W^{(\beta)}_{(n+k)}=0 \;\; (\mbox{respectively, } \;
W^{(\beta)}_{[n+k]}=0) \;\;
\mbox{ for }\;k\in {\mathbb Z}\;\mbox{
sufficiently
negative};\label{lower}\\
&&\dim W^{(\beta)}_{(n)} <\infty \;\; (\mbox{respectively, } \;
\dim W^{(\beta)}_{[n]} <\infty);\nn
&&v_l W^{(\beta)} \subset W^{(\alpha+\beta)}\;\;\mbox{ for any
}\;v\in V^{(\alpha)},\;l\in {\mathbb Z}.\nno
\end{eqnarray}
A strongly $\tilde{A}$-graded (generalized) $V$-module $W$ is said to be
{\it lower bounded} if instead of (\ref{lower}), it satisfies the stronger condition
that for any $\beta \in \tilde{A}$,
\[
W_{(n)}^{(\beta)} = 0\;\; (\mbox{respectively,}\; W_{[n]}^{(\beta)} = 0)
\;\;\mbox{for} \;\; n \in \C\;\; \mbox{and}\;\; \mathfrak{R}(n)\;\;\mbox {sufficiently negative}.
\]
}
\end{defn}

This completes the definition of the notion of strongly
$\tilde{A}$-graded generalized module for a strongly $A$-graded conformal vertex
algebra.

\begin{rema}
{\rm In the strongly graded case, subalgebras (submodules) are vertex subalgebras (submodules) that are strongly graded; algebra and module homomorphisms are of
course understood to preserve the grading by $A$ or $\tilde{A}$.}
\end{rema}

\begin{defn}\label{DH}{\rm Let $V$ be a strongly $A$-graded conformal vertex
algebra. The subspaces $V_{(n)}^{(\alpha)}$ for $n \in \mathbb{Z}$, $\alpha
\in A$ are called the {\it doubly
homogeneous subspaces} of $V$. The elements in $V_{(n)}^{(\alpha)}$
are called {\it doubly homogeneous} elements. Similar definitions
can be used for $W^{(\beta)}_{(n)}$ (respectively,  $W^{(\beta)}_{[n]}$) in the strongly graded (generalized) module $W$.}
\end{defn}

\begin{nota}\label{A-wt}
{\rm Let $v$ be a doubly homogeneous element of $V$. Let wt $v_n$, $n \in \mathbb{Z}$, refer to the weight of $v_n$ as an operator acting on $W$, and let $A$-wt $v_n$ refer to the $A$-weight of $v_n$ on $W$. Similarly, let $w$ be a doubly homogeneous element of $W$. We use wt $w$ to denote the weight of $w$ and $\tilde{A}$-wt $w$ to denote the $\tilde{A}$-grading of $w$.}
\end{nota}

\begin{lemma}\label{b}
Let $v \in V^{(\alpha)}_{(n)}$, for $n \in \mathbb{Z}$,
$\alpha \in A$. Then for $m \in \mathbb{Z}$, {\rm wt}$\ v_m$ = $n - m -1$ and $A$-{\rm wt}$\ v_m = \alpha$.
\end{lemma}
\pf The first equation is standard from the theory of graded conformal vertex algebras
and the second follows easily from the definitions. \epfv

With the strong gradedness condition on a (generalized) module, we can now define the corresponding notion of contragredient module.
\begin{defn}
{\rm
Let $W = \coprod_{\beta \in \tilde{A}, n \in \C} W_{[n]}^{(\beta)}$ be a strongly $\tilde{A}$-graded generalized
module for a strongly $A$-graded conformal vertex algebra. For each $\beta \in \tilde{A}$ and $n \in \C$, let us
identify $(W_{[n]}^{(\beta)})^{*}$ with the subspace of $W^{*}$ consisting of the linear function on $W$ vanishing
on each $W_{[n]}^{(\gamma)}$ with $\gamma \neq \beta$ or $m \neq n$. We define $W'$ to be the $(\tilde{A} \times \C)$-graded
vector subspaces of $W^{*}$ given by
\[
W' =  \coprod_{\beta \in \tilde{A}, n \in \C} (W')_{[n]}^{(\beta)}, \;\; \mbox{where}\;\;  (W')_{[n]}^{(\beta)} =  (W_{[n]}^{(-\beta)})^{*}.
\]
}
\end{defn}

The {\it adjoint vertex operators} $Y'(v, z)\; (v \in V)$ on $W'$ is defined in the same way as vertex operator algebra in Section 5.2 in \cite{FHL}:
\begin{equation}\label{adjoint}
\langle Y'(v, z)w', w\rangle = \langle w', Y(e^{zL(1)}(-z^{-2})^{L(0)}v, z^{-1})w\rangle
\end{equation}
For $w' \in W', w\in W$. The pair $(W', Y')$ carries a strongly graded module structure as follows:

\begin{prop}
Let $\tilde{A}$ be an abelian group containing $A$ as a subgroup and $V$ a strongly $A$-graded conformal vertex algebra. Let $(W, Y)$ be a strongly $\tilde{A}$-graded $V$-module (respectively, generalized $V$-module). Then the pair $(W', Y')$ carries a strongly $\tilde{A}$-graded $V$-module (respectively, generalized $V$-module) structure. If $W$ is lower bounded, so is $W'$.
\end{prop}

\begin{defn}\label{contragredient}
{\rm The pair $(W', Y')$ is called the {\it contragredient module} of $(W, Y)$.}
\end{defn}

\section{$C_1$-cofiniteness condition}
In this Section, we will let $V$ denote a strongly $A$-graded conformal vertex algebra and let $W$ denote a strongly $\tilde{A}$-graded generalized $V$-module, where $A$, $\tilde{A}$ are abelian groups and $A$ is an abelian subgroup of $\tilde{A}$.

In the following definition, we generalize the $C_1$-cofiniteness condition for generalized modules for a vertex operator algebra to a $C_1$-cofiniteness condition for strongly $\tilde{A}$-graded generalized modules for a strongly graded conformal vertex algebra.

\begin{defn}
{\rm Let $C_1(W)$ be the subspace of $W$ spanned by elements of the form $u_{-1}w$ for
\[
u \in V_{+} = \coprod_{n > 0}V_{(n)}
\]
and $w \in W$. The $\tilde{A}$-grading on $W$ induces an $\tilde{A}$-grading on $W/C_1(W)$:
\[
W/C_1(W) = \coprod_{\beta \in \tilde{A}}(W/C_1(W))^{(\beta)},
\]
where
\[
(W/C_1(W))^{(\beta)} = W^{(\beta)}/(C_1(W))^{(\beta)}
\]
for $\beta \in \tilde{A}$.
If dim $(W/C_1(W))^{(\beta)} < \infty$ for $\beta \in \tilde{A}$, we say that $W$ is $C_1$-{\it cofinite with respect to $\tilde{A}$} or $W$ satisfies the $C_1$-{\it cofiniteness condition with respect to $\tilde{A}$}. }
\end{defn}

\begin{nota}{\rm In this paper, unless otherwise stated, we will fixed the grading group for the strongly graded modules. If a strongly graded module is $C_1$-cofinite with respect to $\tilde{A}$, we say it is a $C_1$-cofinite strongly graded module in abbreviation.}
\end{nota}

We associate each $w \in W^{(\beta)}$ with a set of partitions of $\beta$ in $\tilde{A}$, which we denote by $P_w(\beta)$, in the following way:
If $w \notin C_1(W)$, then define $P_w(\beta) = \{\beta\}$;
If $w \in C_1(W)$, then $w$ is a linear span of elements of the form $(u_i)_{-1}w_i'$ for $i = 1, \dots, n$, where $u_i \in V$ is $A$-homogeneous with $A$-$\wt u_i = \alpha_i$ and $w_i' \in W^{(\beta - \alpha_i)}$. We define $P_w(\beta)$ to be the set of partitions $\beta$ in $\tilde{A}$ of the form
\[
\{\alpha_i\}\sqcup P_{w_i'}(\beta - \alpha_i).
\]
We call the set $P_w(\beta)$ an {\it $\tilde{A}$-pattern} of $w$. Note that $P_w(\beta)$ is a finite set and each partition in $P_w(\beta)$ is also a well-defined finite subset of $\tilde{A}$ since $\wt (u_i)_{-1} > 0$. For each $w$, its $\tilde{A}$-pattern may not be unique, but we can always fix an $\tilde{A}$-pattern of $w$.

Without loss of generality, we will assume $P_w(\beta)$ consists of only one element, that is, $P_w(\beta)$ is a partition of $\beta$ in $\tilde{A}$. Also, zeros in $P_w(\beta)$ will not affect our main result, we will still use $P_w(\beta)$ to denote the set of nonzero elements in $P_w(\beta)$.

Let $P = \{\beta_1, \dots, \beta_n\}$ be a partition of $\beta$ in $\tilde{A}$ (note that $\beta_i$ can be the same as $\beta_j$ for $1 \leq i \neq j \leq n$). We define $\mathcal{S}(P)$ to be the set of partitions of $\beta$ consisting of the partitions $\{\gamma_1, \dots, \gamma_m\}$, where $\gamma_i$ is a sum of $\beta_i$'s in $P$ such that the summands of $\gamma_i$ exhaust all the elements in $P$. Note that $P \in \mathcal{S}(P)$ and $\mathcal{S}(P)$ is a finite set if $P$ is a finite set.

\begin{rema}{\rm The set $\mathcal{S}(P)$ has one-to-one correspondence to the set partitions of the set $P$. The cardinality is called the {\it Bell number}.}
\end{rema}

In the remaining context of this paper, we shall study a category of strongly graded generalized $V$-modules satisfying the following condition:
\begin{defn}\label{basicassum}{\rm
For a $C_1$-cofinite strongly graded generalized $V$-module $W$, we say that vertex operators preserve the $\tilde{A}$-patterns of elements in $W$ if for $u \in V^{(\alpha)}$, $w \in W^{(\beta)}$ with $\tilde{A}$-pattern $P_{w}(\beta)$ and $k \in \Z$, the $\tilde{A}$-pattern of $u_kw$ satisfies that
\begin{eqnarray*}
P_{u_kw}(\alpha + \beta) \in
\mathcal{S}(\{\alpha\} \sqcup P_{w}(\beta)).
\end{eqnarray*}}
\end{defn}

\begin{exam}\label{Vitself}{\rm
Let $V$ be a strongly $A$-graded vertex algebra with only nonnegative weights. Then $V$ is $C_1$-cofinite with respect to $A$ as a strongly $A$-graded $V$-module itself since for any $u \in V$, $u = u_{-1}{\bf 1}$.

For any $u \in V^{(\alpha)}$, define the $A$-pattern of $u$ by $P_{u}(\alpha) = \{\alpha\}$. Thus it is easy to see that for $u \in V^{(\alpha)}, v \in V^{(\beta)}$ and $k \in \Z$,
\[
P_{u_kv}(\alpha + \beta) = \{\alpha + \beta\} \in
\mathcal{S}(\{\alpha, \beta\}) = \mathcal{S}(\{\alpha\} \sqcup P_{v}(\beta)).
\]
Hence the vertex operators preserve the $A$-patterns of elements in $V$ defined in the above way.}
\end{exam}

\begin{exam}{\rm
Let $V_L$ be the conformal vertex algebra associated with a nondegenerate even lattice $L$ and let $W$ be a strongly $M$-graded $V_L$-module for a sublattice $M$ of $L^{\circ}$ containing $L$ (see Example \ref{lattice vertex algebra}). Then $W$ is $C_1$-cofinite with respect to $M$.

For the doubly homogeneous elements $w = h_1(-n_1)\cdots h_k(-n_k)\otimes \iota(a) \in W^{(\beta)}$, we can define the $M$-pattern of $w$ by
\[
P_{w}(\beta) = \{\beta\}.
\]

Let $u \in V^{(\alpha)}$. Since $u_kw$ will be a linear combination of elements of the form
\[
g_1(-i_1)\cdots g_m(-i_m)\otimes \iota(b),
\]
where $\bar{b} = \alpha + \beta$, we have
\[
P_{u_kw}(\alpha + \beta) = \{\alpha + \beta\} \in \mathcal{S}(\{\alpha, \beta\}) = \mathcal{S}(\{\alpha\} \sqcup P_{w}(\beta)).
\]
Thus the vertex operators preserve the $M$-patterns of elements in $W$ defined in the above way.}
\end{exam}

\setcounter{equation}{0}
\section{Logarithmic intertwining operators}
 Throughout this paper, we shall use $x, x_0, x_1, x_2, \dots$ to denote commuting formal variables and $z, z_0, z_1, z_2, \dots$ to denote complex variables or complex numbers. We first recall the following definitions.

\begin{defn}\label{log:def}{\rm
Let $(W_1,Y_1)$, $(W_2,Y_2)$ and $(W_3,Y_3)$ be generalized modules
for a conformal vertex algebra $V$. A {\em logarithmic intertwining
operator of type ${W_3\choose W_1\,W_2}$} is a linear map
\begin{equation}\label{log:map0}
{\cal Y}(\cdot, x)\cdot: W_1\otimes W_2\to W_3[\log x]\{x\},
\end{equation}
or equivalently,
\begin{equation}\label{log:map}
w_{(1)}\otimes w_{(2)}\mapsto{\cal Y}(w_{(1)},x)w_{(2)}=\sum_{n\in
{\mathbb C}}\sum_{k\in {\mathbb N}}{w_{(1)}}_{n;\,k}^{\cal
Y}w_{(2)}x^{-n-1}(\log x)^k\in W_3[\log x]\{x\}
\end{equation}
for all $w_{(1)}\in W_1$ and $w_{(2)}\in W_2$, such that the
following conditions are satisfied: the {\em lower truncation
condition}: for any $w_{(1)}\in W_1$, $w_{(2)}\in W_2$ and $n\in
{\mathbb C}$,
\begin{equation}\label{log:ltc}
{w_{(1)}}_{n+m;\,k}^{\cal Y}w_{(2)}=0\;\;\mbox{ for }\;m\in {\mathbb
N} \;\mbox{ sufficiently large,\, independently of}\;k;
\end{equation}
the {\em Jacobi identity}:
\begin{eqnarray}\label{log:jacobi}
\lefteqn{\dps x^{-1}_0\delta \bigg( {x_1-x_2\over x_0}\bigg)
Y_3(v,x_1){\cal Y}(w_{(1)},x_2)w_{(2)}}\nno\\
&&\hspace{2em}- x^{-1}_0\delta \bigg( {x_2-x_1\over -x_0}\bigg)
{\cal Y}(w_{(1)},x_2)Y_2(v,x_1)w_{(2)}\nno \\
&&{\dps = x^{-1}_2\delta \bigg( {x_1-x_0\over x_2}\bigg) {\cal
Y}(Y_1(v,x_0)w_{(1)},x_2) w_{(2)}}
\end{eqnarray}
for $v\in V$, $w_{(1)}\in W_1$ and $w_{(2)}\in W_2$ (note that the
first term on the left-hand side is meaningful because of
(\ref{log:ltc})); the {\em $L(-1)$-derivative property:} for any
$w_{(1)}\in W_1$,
\begin{equation}\label{log:L(-1)dev}
{\cal Y}(L(-1)w_{(1)},x)=\frac d{dx}{\cal Y}(w_{(1)},x).
\end{equation}}
\end{defn}

\begin{defn}\label{slog:def}{\rm In the setting of Definition \ref{log:def}, suppose in addition that $V$ and $W_1$, $W_2$ and $W_3$ are strongly graded. A logarithmic intertwining operator $\cal{Y}$ as in Definition \ref{log:def} is a {\it grading-compatible logarithmic intertwining operator} if for $\beta, \gamma \in \tilde{A}$ and $w_1 \in W_1^{(\beta)}$, $w_2 \in W_2^{(\gamma)}$, $n \in \C$ and $k \in \N$, we have
\[
(w_1)_{n; k}w_2 \in W_3^{(\beta + \gamma)}.
\]
}
\end{defn}

\begin{defn}
{\rm In the setting of Definition \ref{slog:def}, the grading-compatible logarithmic intertwining operators of a fixed type ${W_3\choose W_1\,W_2}$ form a vector space, which we denote by $\mathcal{V}_{W_1 W_2}^{W_3}$. We call the dimension of $\mathcal{V}_{W_1 W_2}^{W_3}$ the {\it fusion rule} for $W_1$, $W_2$ and $W_3$ and denote it by $N_{W_1 W_2}^{W_3}$.}
\end{defn}

\setcounter{equation}{0}
\section{Differential equations}
In the rest of this paper, we assume that $V$ is a strongly $A$-graded vertex algebra, every strongly $\tilde{A}$-graded (generalized) $V$-module is $\R$-graded and satisfies the $C_1$-cofiniteness condition. We also assume the $\tilde{A}$-patterns of elements in the strongly $\tilde{A}$-graded modules are preserved by the vertex operators.

Let $W_i$ be strongly $\tilde{A}$-graded generalized $V$-modules and $\bar{w}_i$ be fixed elements in $W_i^{(\beta_i)}$ for $\beta_i \in \tilde{A}$. Let $P_{\bar{w}_i}(\beta_i)$ be the $\tilde{A}$-pattern of $\bar{w}_i$ for $i = 0, 1, 2, 3$. Set $\beta = \sum_{i=0}^3 \beta_i$ and
\[
P(\bar{w}_0, \bar{w}_1, \bar{w}_2, \bar{w}_3) = \bigsqcup_{i=0}^4 P_{\bar{w}_i}^{\circ}(\beta_i).
\]
Obviously, $P(\bar{w}_0, \bar{w}_1, \bar{w}_2, \bar{w}_3)$ is a partition of $\beta$ in $\tilde{A}$. Form a set
\[
I(\bar{w}_0, \bar{w}_1, \bar{w}_2, \bar{w}_3) = \{(\tilde{\beta}_0, \tilde{\beta}_1, \tilde{\beta}_2, \tilde{\beta}_3)\;|\;\{\tilde{\beta}_0, \tilde{\beta}_1, \tilde{\beta}_2, \tilde{\beta}_3\} \in \mathcal{S}(P(\bar{w}_0, \bar{w}_1, \bar{w}_2, \bar{w}_3))\}.
\]
For simplicity, we will write $P$ for $P(\bar{w}_0, \bar{w}_1, \bar{w}_2, \bar{w}_3)$ and $I$ for $I(\bar{w}_0, \bar{w}_1, \bar{w}_2, \bar{w}_3)$, respectively. The set $I$ is a finite set since $P$ is a finite set.


Let $R = \C[z_1^{\pm 1}, z_2^{\pm 1}, (z_1 - z_2)^{-1}]$. Set
\[
\tilde{T} = \coprod_{(\tilde{\beta}_0, \tilde{\beta}_1, \tilde{\beta}_2, \tilde{\beta}_3) \in I} R \otimes W_0^{(\tilde{\beta}_0)} \otimes W_1^{(\tilde{\beta}_1)} \otimes W_2^{(\tilde{\beta}_2)} \otimes W_3^{(\tilde{\beta}_3)}.
\]
Then $\tilde{T}$ has a natural $R$-module structure.

Let $T$ be an $R$-submodule of $\tilde{T}$ generated by the elements:
\[
\{w_0 \otimes w_1 \otimes w_2 \otimes w_3 \in \tilde{T}| \sqcup _{i=0}^4 P_{w_i}(\tilde{\beta_i}) \in \mathcal{S}(P)\}.
\]

For simplicity, we shall omit one tensor symbol to write $f(z_1, z_2) \otimes w_0 \otimes w_1 \otimes w_2 \otimes w_3$ as $f(z_1, z_2)w_0 \otimes w_1 \otimes w_2 \otimes w_3$ in $\tilde{T}$. For a strongly $\tilde{A}$-graded generalized $V$-module $W$, let $(W', Y')$ be the contragredient module of $W$ (recall definition \ref{contragredient}). In particular, for $u \in V$ and $n \in \Z$, we have the operators $u_n$ on $W'$. Let $u_n^*: W \rightarrow W$ be the adjoint of $u_n: W' \rightarrow W'$. From definition \ref{contragredient} and equation (\ref{adjoint}), we have
\[
\wt u_n^{*} = -\wt u_n
\]
and
\[
\tilde{A}\text{-}\wt u_n^{*} = \tilde{A}\text{-}\wt u_n = \tilde{A}\text{-}\wt u.
\]

\begin{defn}\label{defideal}
{\rm
Let $\alpha$ be a sum of elements in $P$ and $\tilde{\beta}_i \in \tilde{A}$ such that $(\tilde{\beta}_0, \alpha + \tilde{\beta}_1, \tilde{\beta}_2, \tilde{\beta}_3) \in I$. For $u \in V_{+}^{(\alpha)}$, $w_i \in W_i^{(\tilde{\beta}_i)}$ with $\tilde{A}$-patterns of $P_{w_i}(\tilde{\beta}_i)$ of $w_i$ $(i = 0, 1, 2, 3)$ such that
\[
\{\alpha\}\sqcup (\sqcup_{i=0}^3P_{w_i}(\tilde{\beta}_i)) \in \mathcal{S}(P),
\]
we define $J$ to be the submodule of $\tilde{T}$ generated by elements of the form
\begin{eqnarray*}
&&\mathcal{A}(u, w_0, w_1, w_2, w_3)\nno \\
&=& \sum_{k \geq 0}\left(\begin{array}{c}-1\\ k \end{array}\right)(-z_1)^ku_{-1-k}^{*}w_0 \otimes w_1 \otimes w_2 \otimes w_3 - w_0 \otimes u_{-1}w_1 \otimes w_2 \otimes w_3 \nno\\
&& -\sum_{k \geq 0}\left(\begin{array}{c}-1\\ k \end{array}\right)(-(z_1 - z_2))^{-1-k}w_0 \otimes w_1 \otimes u_kw_2 \otimes w_3\nno\\
&& -\sum_{k \geq 0}\left(\begin{array}{c}-1\\ k \end{array}\right)(-z_1)^{-1-k}w_0 \otimes w_1 \otimes w_2 \otimes u_kw_3,\nno\\
&&\mathcal{B}(u, w_0, w_1, w_2, w_3)\nno \\
&=& \sum_{k \geq 0}\left(\begin{array}{c}-1\\ k \end{array}\right)(-z_2)^ku_{-1-k}^{*}w_0 \otimes w_1 \otimes w_2 \otimes w_3 \nno\\
&& -\sum_{k \geq 0}\left(\begin{array}{c}-1\\ k \end{array}\right)(-(z_1 - z_2))^{-1-k}w_0 \otimes u_kw_1 \otimes w_2 \otimes w_3  - w_0 \otimes w_1 \otimes u_{-1}w_2 \otimes w_3\nno\\
&& -\sum_{k \geq 0}\left(\begin{array}{c}-1\\ k \end{array}\right)(-z_2)^{-1-k}w_0 \otimes w_1 \otimes w_2 \otimes u_kw_3,\nno\\
\end{eqnarray*}
\begin{eqnarray*}
&&\mathcal{C}(u, w_0, w_1, w_2, w_3)\nno \\
&=& u_{-1}^{*}w_0 \otimes w_1 \otimes w_2 \otimes w_3 - \sum_{k \geq 0}\left(\begin{array}{c}-1\\ k \end{array}\right)z_1^{-1-k}w_0 \otimes u_kw_1 \otimes w_2 \otimes w_3 \nno\\
&& -\sum_{k \geq 0}\left(\begin{array}{c}-1\\ k \end{array}\right)z_2^{-1-k}w_0 \otimes w_1 \otimes u_kw_2 \otimes w_3 - w_0 \otimes w_1 \otimes w_2 \otimes u_{-1}w_3,\nno\\
&&\mathcal{D}(u, w_0, w_1, w_2, w_3)\nno \\
&=& u_{-1}w_0 \otimes w_1 \otimes w_2 \otimes w_3 \nno\\
&& -\sum_{k \geq 0}\left(\begin{array}{c}-1\\ k \end{array}\right)z_1^{k+1}w_0 \otimes e^{z_1^{-1}L(1)}(-z_1^2)^{L(0)}u_k(-z_1^{-2})^{L(0)}e^{-z_1^{-1}L(1)}w_1 \otimes w_2 \otimes w_3\nno\\
&& -\sum_{k \geq 0}\left(\begin{array}{c}-1\\ k \end{array}\right)z_2^{k+1}w_0 \otimes w_1 \otimes e^{z_2^{-1}L(1)}(-z_2^2)^{L(0)}u_k(-z_2^{-2})^{L(0)}e^{-z_2^{-1}L(1)}w_2 \otimes w_3\nno\\
&& -w_0 \otimes w_1 \otimes w_2 \otimes u_{-1}^{*}w_3.
\end{eqnarray*}}
\end{defn}

\begin{lemma}\label{T}
Let $J$ be the $R$-submodule of $\tilde{T}$ defined in definition \ref{defideal}. Then $J$ is an $R$-submodule of $T$.
\end{lemma}
\pf We show only the case that $\mathcal{A}(u,w_0,w_1,w_2,w_3)$ is in $T$; the other cases are similar. Since the $\tilde{A}$-patterns of elements in $W_0$ are preserved by the vertex operators, we have
\[
P_{u_{-1-k}^{*}w_0}(\alpha + \tilde{\beta_0}) \in \mathcal{S}(\{\alpha\} \sqcup P_{w_0}(\tilde{\beta_0})).
\]
Hence
\[
P_{u_{-1-k}^{*}w_0}(\alpha + \tilde{\beta_0}) \sqcup P_{w_1}(\tilde{\beta}_1)\sqcup P_{w_2}(\tilde{\beta}_2)\sqcup P_{w_3}(\tilde{\beta}_3) \in \mathcal{S}(\{\alpha\}\sqcup (\sqcup_{i=0}^3P_{w_i}(\tilde{\beta}_i)))  \subset \mathcal{S}(P).
\]
Therefore by the definition of $T$, the first expression in $\mathcal{A}(u,w_0,w_1,w_2,w_3)$ lies in $T$. It is similar to show that the other three expressions are also in $T$.\epfv

For $r \in R$, we can define the $R$-submodules $T_{(r)}$, $F_r(T)$ and $F_r(J)$ as in \cite{H}. Note that $F_r(T)$ is a finitely generated $R$-module since $I$ is a finite set.

\begin{prop}\label{main proposition}
There exists $N \in \Z$ such that for any $r \in \R$, $F_r(T) \subset F_r(J) + F_N(T)$. In particular, $T = J + F_{N}(T)$.
\end{prop}
\pf Since $W_i$ is $C_1$-cofinite with respect to $\tilde{A}$, there exists $N_i \in \Z$ such that for $n_i \geq N_i$, $(W_i)_{(n_i)}^{(\beta_i)} \subset (C_1(M_i))^{(\beta_i)}$ for $\beta_i \in \tilde{A}$ ($i = 0, 1, 2, 3$). Let $N = \sum_{i=0}^4 N_i$.

Let $w_0 \otimes w_1 \otimes w_2 \otimes w_3$ be a doubly homogeneous element in $T$. Then we have $w_i \in W_i^{(\tilde{\beta}_i)}$ with $\tilde{A}$-patterns $P_{w_i}(\tilde{\beta}_i)$ such that $\sqcup_{i=0}^3 P_{w_i}(\tilde{\beta}_i) \in \mathcal{S}(P)$. It suffices to show that $w_0 \otimes w_1 \otimes w_2 \otimes w_3 \in  F_N(T) + J$.

If $\wt w_0 \otimes w_1 \otimes w_2 \otimes w_3 \leq N$, the claim is obvious. Now we assume
\[
\wt w_0 \otimes w_1 \otimes w_2 \otimes w_3 > N.
\]
In this case, there exists some $i \in \{0,1,2,3\}$ such that $w_i \in (C_1(W_i))^{(\tilde{\beta}_i)}$. We will discuss the case that $w_1 \in (C_1(W_1))^{(\tilde{\beta}_1)}$, the other cases are completely similar.

Without loss of generality, we assume that $w_1 = u_{-1}w_1'$, where $u$ and $w_1'$ are doubly homogeneous such that $\tilde{A}$-$\wt u = \alpha \in P_{w_1}(\tilde{\beta}_1)$ and $P_{w_1}(\tilde{\beta}_1) = \{\alpha\} \sqcup P_{w_1'}(\tilde{\beta}_1-\alpha)$. By the definition of $\mathcal{A}(u, w_0, w_1', w_2, w_3)$,
\begin{eqnarray*}
&&w_0 \otimes w_1 \otimes w_2 \otimes w_3\nno \\
&=& \sum_{k \geq 0}\left(\begin{array}{c}-1\\ k \end{array}\right)(-z_1)^ku_{-1-k}^{*}w_0 \otimes w_1' \otimes w_2 \otimes w_3 \nno\\
&& -\sum_{k \geq 0}\left(\begin{array}{c}-1\\ k \end{array}\right)(-(z_1 - z_2))^{-1-k}w_0 \otimes w_1' \otimes u_kw_2 \otimes w_3\nno\\
&& -\sum_{k \geq 0}\left(\begin{array}{c}-1\\ k \end{array}\right)(-z_1)^{-1-k}w_0 \otimes w_1' \otimes w_2 \otimes u_kw_3 - \mathcal{A}(u, w_0, w_1', w_2, w_3).
\end{eqnarray*}
Since the $\tilde{A}$-patterns of elements in $W_i$ are preserved by the vertex operators, we have
\[
P_{u_kw_i}(\alpha + \tilde{\beta}_{i}) \in \mathcal{S}(\{\alpha\} \sqcup P_{w_i}(\tilde{\beta}_i)),\;\;
P_{u_{-1-k}^{*}w_0}(\alpha + \tilde{\beta}_{0}) \in \mathcal{S}(\{\alpha\} \sqcup P_{w_i}(\tilde{\beta}_0)).
\]
Using the same argument as in the proof of lemma \ref{T}, we see that the elements $u_{-1-k}^{*}w_0 \otimes w_1' \otimes w_2 \otimes w_3$, $w_0 \otimes w_1' \otimes u_kw_2 \otimes w_3$ and $w_0 \otimes w_1' \otimes w_2 \otimes u_kw_3$ are all in $T$. Note that these elements all have weights strictly less than $w_0 \otimes w_1 \otimes w_2 \otimes w_3$, thus the element $w_0 \otimes w_1 \otimes w_2 \otimes w_3$ can be written as a sum of elements of $J$ and elements of $T$ with lower weights. By induction on the weight of $w_0 \otimes w_1 \otimes w_2 \otimes w_3$, we know that $w_0 \otimes w_1 \otimes w_2 \otimes w_3$ can be written as a sum of elements of $J$ and an element of $F_N(T)$. \epfv

For an element $\mathcal{W} \in T$, we shall use $[\mathcal{W}]$ to denote the equivalence class in $T/J$ containing $\mathcal{W}$. We have the following theorem:

\begin{thm}
Let $W_i$ be strongly $\tilde{A}$-graded generalized $V$-modules for $i = 0, 1, 2, 3$. For any $\tilde{A}$-homogeneous elements $w_i \in W_i$ $(i = 0, 1, 2, 3)$, let $M_1$ and $M_2$ be the $R$-submodules of $T/J$ generated by $[w_0 \otimes L(-1)^jw_1 \otimes w_2 \otimes w_3]$, $j \geq 0$, and by $[w_0 \otimes w_1 \otimes L(-1)^jw_2 \otimes w_3]$, $j \geq 0$, respectively. Then $M_1$, $M_2$ are finitely generated. In particular, for any $\tilde{A}$-homogeneous elements $w_i \in W_i$ $(i = 0, 1, 2, 3)$, there exist $a_k(z_1, z_2)$, $b_l(z_1, z_2) \in R$ for $k = 1, \dots, m$ and $l = 1, \dots, n$ such that
\begin{eqnarray}\label{e1}
[w_0 \otimes L(-1)^mw_1 \otimes w_2 \otimes w_3] + a_1(z_1, z_2)[w_0 \otimes L(-1)^{m-1}w_1 \otimes w_2 \otimes w_3]\nno\\
+ \cdots + a_m(z_1, z_2)[w_0 \otimes w_1 \otimes w_2 \otimes w_3] = 0,
\end{eqnarray}

\begin{eqnarray}\label{e2}
[w_0 \otimes w_1 \otimes L(-1)^nw_2 \otimes w_3] + b_1(z_1, z_2)[w_0 \otimes w_1 \otimes L(-1)^{n-1}w_2 \otimes w_3]\nno\\
+ \cdots + b_n(z_1, z_2)[w_0 \otimes w_1 \otimes w_2 \otimes w_3] = 0.
\end{eqnarray}
\end{thm}
\pf We construct $R$-module $T$ and $R$-submodule $J \subset T$ associated to for each quadruple $w_i \in W_i^{(\beta_i)}$ ($i = 0, 1, 2, 3$). By proposition \ref{main proposition}, $T/J$ is finitely generated. Since $R$ is a Noetherian ring, any $R$-submodule of the finitely generated $R$-module $T/J$ is also finitely generated. In particular, $M_1$ and $M_2$ are finitely generated. The second conclusion follows immediately. \epfv

Now we establish the existence of systems of differential equations:
\begin{thm}\label{main theorem}
Let $W_i$ for $i = 0, 1, 2, 3$ be strongly $\tilde{A}$-graded generalized $V$-modules satisfying the $C_1$-cofiniteness condition. Suppose that the $\tilde{A}$-patterns of elements in $W_i$ are preserved by the vertex operators. Then for any $\tilde{A}$-homogeneous elements $w_i \in W_i$ ($i = 0, 1, 2, 3$), there exist
\[
a_k(z_1, z_2), b_l(z_1, z_2) \in \C[z_1^{\pm}, z_2^{\pm}, (z_1 - z_2)^{-1}]
\]
for $k = 1, \dots, m$ and $l = 1, \dots, n$ such that for any strongly graded $V$-module $W_4, W_5$ and $W_6$, any intertwining operators $\mathcal{Y}_1, \mathcal{Y}_2, \mathcal{Y}_3, \mathcal{Y}_4, \mathcal{Y}_5$ and $\mathcal{Y}_6$ of types ${W_0'\choose W_1\,W_4}$, ${W_4\choose W_2\,W_3}$, ${W_5\choose W_1\,W_2}$, ${W_0'\choose W_5\,W_3}$, ${W_0'\choose W_2\,W_6}$ and ${W_6\choose W_1\,W_3}$, respectively, the series
\begin{equation}\label{e3}
\langle w_0, \mathcal{Y}_1(w_1, z_1)\mathcal{Y}_2(w_2, z_2)w_3\rangle,
\end{equation}
\begin{equation}\label{e4}
\langle w_0, \mathcal{Y}_4(\mathcal{Y}_3(w_1, z_1-z_2)w_2, z_2)w_3\rangle
\end{equation}
and
\begin{equation}\label{e5}
\langle w_0, \mathcal{Y}_5(w_2, z_2)\mathcal{Y}_6(w_1, z_1)w_3\rangle,
\end{equation}
satisfy the expansions of the system of differential equations
\begin{equation}\label{eq7}
\frac{\partial^m \varphi}{\partial z_1^m} + a_1(z_1, z_2)\frac{\partial^{m-1} \varphi}{\partial z_1^{m-1}} + \cdots + a_m(z_1, z_2)\varphi = 0,
\end{equation}
\begin{equation}\label{eq8}
\frac{\partial^n \varphi}{\partial z_2^n} + b_1(z_1, z_2)\frac{\partial^{n-1} \varphi}{\partial z_2^{n-1}} + \cdots + b_n(z_1, z_2)\varphi = 0
\end{equation}
in the region $|z_1| > |z_2| > 0$, $|z_2|>|z_1-z_2|>0$ and $|z_2|>|z_1|>0$, respectively.
\end{thm}

\pf The proof is similar to the proof of Theorem 1.4 in \cite{H}. We sketch the proof as follows:

Let $\Delta = {\rm wt}\ w_0 - {\rm wt}\ w_1- {\rm wt}\ w_2 - {\rm wt}\ w_3$. Let $\C(\{x\})$ be the space of all series of the form $\sum_{n \in \R}a_nx^n$ for $n \in \R$ such that $a_n = 0$ when the real part of $n$ is sufficiently negative.

Consider the map
\begin{eqnarray*}
\phi_{\mathcal{Y}_1, \mathcal{Y}_2}: T \longrightarrow z_1^{\Delta}\C(\{z_2/z_1\})[z_1^{\pm 1}, z_2^{\pm 1}]
\end{eqnarray*}
defined by
\begin{eqnarray*}
& \phi_{\mathcal{Y}_1, \mathcal{Y}_2}(f(z_1, z_2)w_0 \otimes w_1 \otimes w_2 \otimes w_3) \\
& = \iota_{|z_1| > |z_2| > 0}(f(z_1, z_2))\langle w_0, \mathcal{Y}_1(w_1, z_1)\mathcal{Y}_2(w_2, z_2)w_3\rangle,
\end{eqnarray*}
where
\begin{eqnarray*}
\iota_{|z_1| > |z_2| > 0}: R &\longrightarrow & \C[[z_2/z_1]][z_1^{\pm 1}, z_2^{\pm 1}]
\end{eqnarray*}
is the map expanding elements of $R$ as series in the regions $|z_1| > |z_2| > 0$.

Using the Jacobi identity for the logarithmic intertwining operators, we have that $\phi_{\mathcal{Y}_1, \mathcal{Y}_2}(J) = 0$. Thus the map $\phi_{\mathcal{Y}_1, \mathcal{Y}_2}$ induces a map
\begin{eqnarray*}
\bar{\phi}_{\mathcal{Y}_1, \mathcal{Y}_2}: T/J \longrightarrow z_1^{\Delta}\C(\{z_2/z_1\})[z_1^{\pm 1}, z_2^{\pm 1}].
\end{eqnarray*}
Applying $\bar{\phi}_{\mathcal{Y}_1, \mathcal{Y}_2}$ to (\ref{e1}) and (\ref{e2}) and then use the $L(-1)$-derivative property for logarithmic intertwining operators, we see that (\ref{e3}) indeed satisfies the expansions of the system of differential equations in the regions $|z_1| > |z_2| > 0$. Similarly, we can prove that (\ref{e4}) and (\ref{e5}) satisfy the expansions of the system of differential equations in the regions $|z_2| > |z_1-z_2| > 0$ and $|z_2|>|z_1|>0$, respectively. \epfv

The following result can be proved by the same method, so we omit the proof.
\begin{thm}\label{theorem 2}
Let $W_i$ be strongly $\tilde{A}$-graded generalized $V$-modules satisfying the $C_1$-cofiniteness condition for $i = 0, \dots, n + 1$. Suppose that the $\tilde{A}$-patterns of elements in $W_i$ are preserved by the vertex operators. For any generalized V-modules $\widetilde{W_1}, \dots, \widetilde{W_{n-1}}$, let
\[
\mathcal{Y}_1, \mathcal{Y}_2, \dots, \mathcal{Y}_{n-1}, \mathcal{Y}_n
\]
be logarithmic intertwining operators of types
\[
{W_0\choose W_1\,\widetilde{W_1}}, {\widetilde{W_1}\choose W_2\,\widetilde{W_2}}, \dots, {\widetilde{W_{n-2}}\choose W_{n-1}\,\widetilde{W_{n-1}}}, {\widetilde{W_{n-1}}\choose W_n\,W_{n+1}},
\]
respectively. Then for any $\tilde{A}$-homogeneous elements $w_{(0)}' \in W_0'$, $w_{(1)} \in W_1, \dots, w_{(n + 1)} \in W_{n+1}$, there exist
\[
a_{k, l} (z_1, \dots, z_n) \in \C[z_1^{\pm 1}, \dots, z_n^{\pm 1}, (z_1 - z_2)^{-1}, (z_1 - z_3)^{-1}, \dots, (z_{n-1} - z_n)^{-1}]
\]for $k = 1, \dots, m$ and $l = 1, \dots, n$ such that the series
\[
\langle w_{(0)}', \mathcal{Y}_1(w_{(1)}, z_1)\cdots \mathcal{Y}_n(w_{(n)}, z_n)w_{(n+1)}\rangle
\]
satisfies the system of differential equations
\begin{equation}\label{eq9}
\frac{\partial^m\varphi}{\partial z_l^m} + \sum_{k = 1}^m a_{k, l}(z_1, \dots, z_n)\frac{\partial^{m - k}\varphi}{\partial z_l^{m - k}} = 0, \ \ l = 1, \dots, n
\end{equation}
in the region $|z_1| > \cdots > |z_n| > 0$.
\end{thm}

\setcounter{equation}{0}
\section{The regularity of the singular points}
We first recall the definition for {\it regular singular points} for a system of differential equations given in \cite{K}. For the system of differential equations of form (\ref{eq9}), {\it a singular point}
\[
z_0 = (z_0^{(1)}, \dots, z_0^{(n)})
\]
is an isolated singular point of the coefficient matrix
\[
a_{k, l} (z_1, \dots, z_n) \in \C[z_1^{\pm 1}, \dots, z_n^{\pm 1}, (z_1 - z_2)^{-1}, (z_1 - z_3)^{-1}, \dots, (z_{n-1} - z_n)^{-1}]
\]
for $k = 1, \dots, m$ and $l = 1, \dots, n$.
For $s = (s_1, \dots, s_n) \in \Z_{+}^n$, set
\[
|s| = \sum_{i=0}^n s_i
\]
and
\[
(\log (z-z_0))^s = (\log (z_1 - z_0^{(1)}))^{s_1}\cdots (\log (z_n - z_0^{(n)}))^{s_n}.
\]
For $t = (t^{(1)}, \dots, t^{(n)}) \in \C^n$, set
\[
(z-z_0)^t = (z_1-z_0^{(1)})^{t^{(1)}}\cdots (z_n-z_0^{(n)})^{t^{(n)}}.
\]
A singular point $z_0$ for the system of differential equations of form (\ref{eq9}) is {\it regular} if every solution in the punctured disc $(D^{\times})^n$
\[
0< |z_i - z_0^{(i)}| < a_i
\]
with some $a_i \in \R_{+}$ $(i = 1, \dots, n)$ is of the form
\[
\varphi(z) = \sum_{i = 1}^r\sum_{|m| < M}(z-z_0)^{t_i}(\log (z-z_0))^mf_{t_i, m}(z-z_0)
\]
with $M, r \in \Z_{+}$ and each $f_{t_i, m}(z-z_0)$ holomorphic in $(D^{\times})^n$. Theorem B.16 in \cite{K} gives a sufficient condition for a singular point of a system of differential equations to be regular.

As in \cite{H}, for $r \in \R$, we define the $R$-modules $F_r^{(z_1 = z_2)}(R)$, $F_r^{(z_1 = z_2)}(T)$ and $F_r^{(z_1 = z_2)}(\widetilde{T})$, which provide filtration associated to the singular point $z_1 = z_2$ on $R$, $R$-modules $T$ and $\widetilde{T}$, respectively.

Let $F_r^{(z_1 = z_2)}(J) = F_r^{(z_1 = z_2)}(T) \cap J$ for $r \in \R$. We have the following refinement of Proposition \ref{main proposition}:
\begin{prop}\label{l1}
For any $r \in \R$, $F_r^{(z_1 = z_2)}(T) \subset F_r^{(z_1 = z_2)}(J) + F_N(T)$.
\end{prop}
\pf The proof is similar to the proof of proposition \ref{main proposition} except for some slight differences. We discuss elements of the form $w_0 \otimes u_{-1}w_1 \otimes w_2 \otimes w_3$ with weight $s$, where $w_i \in W_i^{(\widetilde{\beta_i})}$ for $i = 0, 1, 2, 3$ and $u \in (V_0)_{+}$. By definition of the element $\mathcal{A}(u, w_0, w_1, w_2, w_3)$ in the $R$-submodule $J$, we have
\begin{eqnarray*}
&& w_0 \otimes u_{-1}w_1 \otimes w_2 \otimes w_3 \nno \\
&& \;\;\;\;\;\;\;\;\;= \sum_{k \geq 0}\left(\begin{array}{c}-1\\ k \end{array}\right)(-z_1)^ku_{-1-k}^{*}w_0 \otimes w_1 \otimes w_2 \otimes w_3 - \mathcal{A}(u, w_0, w_1, w_2, w_3)\nno\\
&& \;\;\;\;\;\;\;\;\;\;\;\;\;\;-\sum_{k \geq 0}\left(\begin{array}{c}-1\\ k \end{array}\right)(-(z_1 - z_2))^{-1-k}w_0 \otimes w_1 \otimes u_kw_2 \otimes w_3\nno\\
&& \;\;\;\;\;\;\;\;\;\;\;\;\;\;-\sum_{k \geq 0}\left(\begin{array}{c}-1\\ k \end{array}\right)(-z_1)^{-1-k}w_0 \otimes w_1 \otimes w_2 \otimes u_kw_3.
\end{eqnarray*}
We know that the elements $u_{-1-k}^{*}w_0 \otimes w_1 \otimes w_2 \otimes w_3$, $w_0 \otimes w_1 \otimes u_kw_2 \otimes w_3$ and $w_0 \otimes w_1 \otimes w_2 \otimes u_kw_3$ for $k \geq 0$ lie in $T$ with weights less than the weight of $w_0 \otimes u_{-1}w_1 \otimes w_2 \otimes w_3$.

By induction assumption, $u_{-1-k}^{*}w_0 \otimes w_1 \otimes w_2 \otimes w_3$, $w_0 \otimes w_1 \otimes w_2 \otimes u_kw_3 \in F_s^{(z_1 = z_2)}(J) + F_N(T)$ and $w_0 \otimes w_1 \otimes u_kw_2 \otimes w_3 \in F_{s - k - 1}^{(z_1 = z_2)}(J) + F_N(T)$. Hence the element $(-(z_1 - z_2))^{-1-k}w_0 \otimes w_1 \otimes u_kw_2 \otimes w_3 \in F_s^{(z_1 = z_2)}(J) + F_N(T)$. Thus in this case, $w_0 \otimes u_{-1}w_1 \otimes w_2 \otimes w_3$ can be written as a sum of an element of $F_s^{(z_1 =z_2)}(J)$ and an element of $F_N(T)$.\epfv

We shall also consider the ring $\C[z_1^{\pm}, z_2^{\pm}]$ and the $\C[z_1^{\pm}, z_2^{\pm}]$-module
\[
T^{(z_1 = z_2)} = \coprod_{(\tilde{\beta}_0, \tilde{\beta}_1, \tilde{\beta}_2, \tilde{\beta}_3) \in I} \C[z_1^{\pm}, z_2^{\pm}] \otimes \{w_0\otimes w_1\otimes w_2 \otimes w_3|w_i \in W_i^{(\tilde{\beta_i})},\; \sqcup_{i=0}^3 P_{w_i}(\tilde{\beta}_i) \in \mathcal{S}(P)\}.
\]
Let $T_{(r)}^{(z_1 = z_2)}$ be the space of elements of $T^{(z_1 = z_2)}$ of weight $r$ for $r \in \R$. Let $F_r(T^{(z_1 = z_2)}) = \coprod_{s \leq r} T_{(s)}^{(z_1 = z_2)}$. These subspaces give a filtration of $T^{(z_1 = z_2)}$ in the following sense:  $F_r(T^{(z_1 = z_2)}) \subset F_s(T^{(z_1 = z_2)})$ for $r \leq s$ and $T^{(z_1 = z_2)} = \coprod_{r \in \R}F_r(T^{(z_1 = z_2)})$.

Let $w_i \in W_i^{(\tilde{\beta}_i)}$, where $(\tilde{\beta}_0, \tilde{\beta}_1, \tilde{\beta}_2, \tilde{\beta}_3) \in I$ and $\sqcup_{i=0}^3P_{w_i}(\tilde{\beta}_i) \in \mathcal{S}(P)$. Then by Proposition \ref{l1},
\[
w_0 \otimes w_1 \otimes w_2 \otimes w_3 = \mathcal{W}_1 + \mathcal{W}_2
\]
where $\mathcal{W}_1 \in F_{\sigma}^{(z_1 = z_2)}(J)$ and $\mathcal{W}_2 \in F_N(T)$. Using the same proof as Lemma 2.2 in \cite{H}, we have the following lemma:

\begin{lemma}\label{l2}
For any $s \in [0, 1)$, there exist $S \in \R$ such that $s + S \in \Z_{+}$ and for any $w_i \in W_i$, $i = 0, 1, 2, 3$, satisfying $\sigma \in s + \Z$, $(z_1 - z_2)^{\sigma + S}\mathcal{W}_2 \in T^{(z_1 = z_2)}$.
\end{lemma}

\begin{thm}\label{regsing}
Let $W_i$, $w_i \in W_i$ for $i = 0, 1, 2, 3, 4$, $\mathcal{Y}_1$ and $\mathcal{Y}_2$ be the same as in Theorem \ref{main theorem}. For any possible singular point of the form $(z_1 = 0, z_2 = 0, z_1 = \infty, z_2 = \infty, z_1 = z_2)$, $z_1^{-1}(z_1 - z_2) = 0$, or $z_2^{-1}(z_1 - z_2) = 0$, there exist
\[
a_k(z_1, z_2), b_l(z_1, z_2) \in \C[z_1^{\pm}, z_2^{\pm}, (z_1 - z_2)^{-1}]
\]
for $k = 1, \dots, m$ and $l = 1, \dots, n$, such that this singular point of the system (\ref{eq7}) and (\ref{eq8}) satisfied by (\ref{e3}) is regular.
\end{thm}

\pf The proof is the same as the proof of Theorem 2.3 in \cite{H} except that we use Proposition \ref{l1} and Lemma \ref{l2} here.\epfv

We can prove the following theorem using the same method, so we omit the proof here.
\begin{thm}\label{theorem 3}
For any set of possible singular points of the system (\ref{eq9}) in Theorem \ref{theorem 2} of the form $z_i = 0$ or $z_i = \infty$ for some $i$ or $z_i = z_j$ for some $i \neq j$, the $a_{k, l}(z_1, \dots, z_n)$ in Theorem \ref{theorem 2} can be chosen for $k = 1, \dots, m$ and $l = 1, \dots, n$ so that these singular points are regular.
\end{thm}

\setcounter{equation}{0}
\section{Convergence and extension property}
In the logarithmic tensor category theory developed in \cite{HLZ1}-\cite{HLZ8}, the convergence and expansion property for the logarithmic intertwining operators are needed in the construction of the associativity isomorphism. In this Section, we will recall the definition of convergence and expansion property for products and iterates of logarithmic intertwining operators and then follow \cite{HLZ7} and \cite{HLZ8} to give sufficient conditions for a category to have these properties.

Throughout this Section, we will let $\mathcal{M}_{sg}$ (respectively, $\mathcal{GM}_{sg}$) denote the category of the strongly $\tilde{A}$-graded (respectively, generalized) $V$-modules. We are going to study the subcategory $\mathcal{C}$ of $\mathcal{M}_{sg}$ (respectively, $\mathcal{GM}_{sg}$) satisfying the following assumptions.

\begin{assum}\label{a3}
{\rm We shall assume the following:
\begin{itemize}
\item $A$ and $\tilde{A}$ are abelian groups satisfying $A \leq \tilde{A}$.
\item $V$ is a strongly $A$-graded conformal vertex algebra and $V$ is an object of $\cal{C}$ as a $V$-module.
\item All (generalized) $V$-modules are lower bounded, satisfy the $C_1$-cofiniteness condition with respect to $\tilde{A}$.
\item The $\tilde{A}$-patterns of elements in the (generalized) modules are preserved by the vertex operators.
\item For any object of $\cal{C}$, the (generalized) weights are real numbers and in addition there exist $K \in \Z$ such that $(L(0) - L(0)_s)^K = 0$ on the generalized module.
\item $\cal{C}$ is closed under images, under the contragredient functor, under taking finite direct sums.
\end{itemize}}
\end{assum}

Given objects $W_1, W_2, W_3, W_4, M_1$ and $M_2$ of the category $\cal{C}$, let $\mathcal{Y}_1, \mathcal{Y}_2, \mathcal{Y}^{1}$ and $\mathcal{Y}^2$ be logarithmic intertwining operators of types ${W_4\choose W_1\,M_1}$, ${M_1\choose W_2\,W_3}$, ${W_4\choose M_2\,W_3}$ and ${M_2\choose W_1\,W_2}$, respectively. We recall the following definitions and theorems from Section $11$ in \cite{HLZ7}:

{\bf Convergence and extension property for products} For any $\beta \in \tilde{A}$, there exists an integer $N_{\beta}$ depending only on $\mathcal{Y}_1$ and $\mathcal{Y}_2$ and $\beta$, and for any doubly homogeneous elements $w_{(1)} \in (W_1)^{(\beta_1)}$ and $w_{(2)} \in (W_2)^{(\beta_2)}$ $(\beta_1, \beta_2 \in \tilde{A})$ and any $w_{(3)} \in W_3$ and $w_{(4)}' \in W_4'$ such that \[\beta_1 + \beta_2 = -\beta,\]there exist $M \in \N$, $r_k, s_k \in \R$, $i_k, j_k \in \N$, $k = 1, \dots, M$, and analytic functions $f_k(z)$ on $|z| < 1$, $k = 1, \dots, M$, satisfying \[\wt w_{(1)} + \wt w_{(2)} + s_k > N_{\beta},\ k = 1, \dots, M,\]such that \[\langle w_{(4)}', \mathcal{Y}_1(w_{(1)},  x_1)\mathcal{Y}_2(w_{(2)}, x_2)w_{(3)} \rangle_{W_4} |_{x_1 = z_1,\ x_2 = z_2}\] is absolutely convergent when $|z_1| > |z_2| > 0$ and can be analytically extended to the multivalued analytic function \[\sum_{k = 1}^M z_2^{r_k}(z_1 - z_2)^{s_k}(\log z_2)^{i_k}(\log (z_1 - z_2))^{j_k}f_k(\frac{z_1 - z_2}{z_2})\](here $\log(z_1 - z_2)$ and $\log z_2$, and in particular, the powers of the variables, mean the multivalued functions, not the particular branch we have been using) in the region $|z_2| > |z_1 - z_2| > 0$.

{\bf Convergence and extension property without logarithms for products} When $i_k = j_k = 0$ for $k = 1, \dots, M$, we call the property above the {\it convergence and extension property without logarithms for products}.

{\bf Convergence and extension property for iterates} For any $\beta \in \tilde{A}$, there exists an integer $\tilde{N_{\beta}}$ depending only on $\mathcal{Y}^1$ and $\mathcal{Y}^2$ and $\beta$, and for any doubly homogeneous elements $w_{(1)} \in (W_1)^{(\beta_1)}$ and $w_{(2)} \in (W_2)^{(\beta_2)}$ $(\beta_1, \beta_2 \in \tilde{A})$ and any $w_{(3)} \in W_3$ and $w_{(4)}' \in W_4'$ such that \[\beta_1 + \beta_2 = -\beta,\]there exist $\tilde{M} \in \N$, $\tilde{r_k}, \tilde{s_k} \in \R$, $\tilde{i_k}, \tilde{j_k} \in \N$, $k = 1, \dots, \tilde{M}$, and analytic functions $\tilde{f_k}(z)$ on $|z| < 1$, $k = 1, \dots, M$, satisfying \[\wt w_{(1)} + \wt w_{(2)} + \tilde{s_k} > \tilde{N_{\beta}},\ k = 1, \dots, \tilde{M},\]such that \[\langle w_{(0)}', \mathcal{Y}_1(\mathcal{Y}_2(w_{(1)}, x_0)w_{(2)}, x_2)w_{(3)} \rangle_{W_4} |_{x_0 = z_1 - z_2,\ x_2 = z_2}\] is absolutely convergent when $|z_2| > |z_1 - z_2| > 0$ and can be analytically extended to the multivalued analytic function \[\sum_{k = 1}^{\tilde{M}} z_1^{\tilde{r_k}}z_2^{\tilde{s_k}}(\log z_1)^{\tilde{i_k}}(\log z_2)^{\tilde{j_k}}\tilde{f_k}(\frac{z_2}{z_1})\](here $\log z_1$ and $\log z_2$, and in particular, the powers of the variables, mean the multivalued functions, not the particular branch we have been using) in the region $|z_1| > |z_2| > 0$.

{\bf Convergence and extension property without logarithmic for iterates} When $i_k = j_k = 0$ for $k = 1, \dots, M$, we call the property above the {\it convergence and extension property without logarithms for iterates}.

If the convergence and extension property (with or without logarithms) for products holds for any objects $W_1, W_2, W_3, W_4$ and $M_1$ of $\cal{C}$ and any logarithmic intertwining operators $\mathcal{Y}_1$ and $\mathcal{Y}_2$ of the types as above, we say that {\it the convergence and extension property for products holds in $\cal{C}$}. We similarly define the meaning of the phrase {\it the convergence and extension property for iterates holds in $\cal{C}$.}

The following theorem generalizes Theorem $11.8$ in \cite{HLZ7} to the strongly graded generalized modules for a strongly graded conformal vertex algebra:
\begin{thm}\label{theorem 1}
Let $V$ be a strongly graded conformal vertex algebra. Then
\begin{itemize}
 \item [1.] The convergence and extension properties for products and iterates hold in $\cal{C}$. If $\cal{C}$ is in $\mathcal{M}_{sg}$ and if every object of $\cal{C}$ is a direct sum of irreducible objects of $\cal{C}$ and there are only finitely many irreducible objects of $\cal{C}$ (up to equivalence), then the convergence and extension properties without logarithms for products and iterates hold in $\cal{C}$.
 \item [2.] For any $n \in \Z_{+}$, any objects $W_1, \dots, W_{n+1}$ and $\widetilde{W_1}, \dots, \widetilde{W_{n-1}}$ of $\cal{C}$, any logarithmic intertwining operators \[\mathcal{Y}_1, \mathcal{Y}_2, \dots, \mathcal{Y}_{n-1}, \mathcal{Y}_n\]of types \[{W_0\choose W_1\,\widetilde{W_1}}, {\widetilde{W_1}\choose W_2\,\widetilde{W_2}}, \dots, {\widetilde{W_{n-2}}\choose W_{n-1}\,\widetilde{W_{n-1}}}, {\widetilde{W_{n-1}}\choose W_n\,W_{n+1}},\]respectively, and any $w_{(0)}' \in W_0'$, $w_{(1)} \in W_1, \dots, W_{(n+1)} \in W_{n+1}$, the series \[\langle w_{(0)}', \mathcal{Y}_1(w_{(1)}, z_1)\cdots \mathcal{Y}_n(w_{(n)}, z_n)w_{(n+1)}\rangle\]is absolutely convergent in the region $|z_1| > \cdots > |z_n| > 0$ and its sum can be analytically extended to a multivalued analytic function on the region given by $z_1 \neq 0$, $i = 1, \dots, n$, $z_i \neq z_j$, $i \neq j$, such that for any set of possible singular points with either $z_i = 0, z_i = \infty$ or $z_i = z_j$ for $i \neq j$, this multivalued analytic function can be expanded near the singularity as a series having the same form as the expansion near the singular points of a solution of a system of differential equations with regular singular points.
\end{itemize}
\end{thm}
\pf The first statement in the first part and the statement in the second part of the theorem follow directly from Theorem \ref{theorem 2} and Theorem \ref{theorem 3} and the theorem of differential equations with regular singular points. The second statement in the first part can be proved using the same method in \cite{H}. \epfv

\setcounter{equation}{0}

\section{Examples}
\begin{exam}\label{vertex operator algebra}{\rm The notion of conformal vertex algebra strongly
graded with respect to the trivial group is exactly the notion of
vertex operator algebra. Let $V$ be a vertex operator algebra,
viewed (equivalently) as a conformal vertex algebra strongly graded
with respect to the trivial group.  Then the $V$-modules that are
strongly graded with respect to the trivial group (in the sense of
Definition \ref{def:dgw}) are exactly the ($\C$-graded) modules for
$V$ as a vertex operator algebra, with the grading restrictions as
follows: For $n \in \C$,
\[
W_{(n+k)}=0 \;\; \mbox { for }\;k\in {\mathbb Z}\;\mbox{
sufficiently negative}
\]
and
\[
\dim W_{(n)} <\infty.
\]
}
\end{exam}

\begin{exam}\label{lattice vertex algebra}{\rm
An important source of examples of strongly graded conformal vertex algebras and modules comes {}from the vertex algebras and modules associated with even lattices. We recall the following construction from \cite{FLM}. Let $L$ be an even lattice, i.e., a finite-rank free abelian group equipped with a nondegenerate symmetric bilinear form $\langle\cdot,\cdot\rangle$, not necessarily
positive definite, such that $\langle \alpha, \alpha\rangle\in 2{\mathbb Z}$
for all $\alpha\in L$. Let $\mathfrak{h} = L\otimes_{\mathbb{Z}}
\mathbb{C}$. Then $\mathfrak{h}$ is a vector space with a nonsingular
bilinear form $\langle \cdot, \cdot \rangle$, extended from $L$. We form a Heisenberg algebra
\[
\widehat{\mathfrak{h}}_{\mathbb{Z}} = \coprod_{n \in \mathbb{Z},\ n \neq 0}
\mathfrak{h} \otimes t^n \oplus \mathbb{C}c.
\]
Let $(\widehat{L},\bar{}\ )$ be a central extension of $L$ by a finite cyclic group $\langle \kappa\;|\;\kappa^s =
 1\rangle$. Fix a primitive $s$th root of unity, say $\omega$, and define
 the faithful character
\[
\chi : \langle \kappa \rangle \rightarrow
\mathbb{C}^{*}
\]
 by the condition
\[
 \chi(\kappa) = \omega.
\]
Denote by $\mathbb{C}_{\chi}$ the one-dimensional space $\mathbb{C}$ viewed as a $\langle \kappa \rangle$-module on which $\langle \kappa \rangle$ acts according to
$\chi$:
\[
\kappa \cdot 1 = \omega,
\]
and denote by $\mathbb{C}\{L\}$ the induced $\widehat{L}$-module
\[
\mathbb{C}\{L\} = {\rm Ind}^{\widehat{L}}_{\langle \kappa \rangle} \mathbb{C}_{\chi} =
\mathbb{C}[\widehat{L}]\otimes_{\mathbb{C}[\langle \kappa \rangle]} \mathbb{C}_{\chi}.
\]
Then
\[
V_L = S(\widehat{\mathfrak{h}}_{\mathbb{Z}}^{-}) \otimes \mathbb{C}\{L\}
\]
has a natural structure of conformal vertex algebra; see \cite{B1} and Chapter 8 of
\cite{FLM}. For $\alpha \in L$, choose an $a \in \widehat{L}$ such that
$\bar{a} = \alpha$. Define
\[
\iota(a) = a \otimes 1 \in \mathbb{C}\{L\}
\]
and
\[
V_L^{(\alpha)} = \mbox{span}\; \{h_1(-n_1)\cdots h_k(-n_k)\otimes \iota(a)\},
\]
where $h_1, \dots, h_k \in \mathfrak{h}$,
$n_1, \dots, n_k > 0$, and where $h(n)$ is the natural operator associated with $h \otimes t^n$
via the $\hat{\mathfrak{h}}_{\Z}$-module structure of $V_L$. Then $V_L$ is equipped with a natural second grading
given by $L$ itself. Also for $n \in \mathbb{Z}$, we have
\[
(V_L)_{(n)}^{(\alpha)} = {\rm span}\ \{h_1(-n_1)\cdots
h_k(-n_k)\otimes \iota(a)|\ \bar{a} = \alpha, \sum_{i = 1}^k n_i + \frac{1}{2}\langle \alpha, \alpha\rangle = n\},
\]
making $V_L$ a strongly $L$-graded conformal vertex algebra in the sense of definition \ref{def:dgv}. When the form $\langle\cdot,\cdot\rangle$ on $L$ is also positive definite, then
$V_L$ is a vertex operator algebra, that is, as in example \ref{vertex operator algebra}, $V_L$ is a strongly
$A$-graded conformal vertex algebra for $A$ the trivial group. In general, a conformal vertex algebra
may be strongly graded for several choices of $A$.

Any sublattice $M$ of the ``dual lattice''
$L^{\circ}$ of $L$ containing $L$ gives rise to a strongly
$M$-graded module for the strongly $L$-graded conformal vertex
algebra (see Chapter 8 of \cite{FLM}; cf. \cite{LL}). In fact, any
irreducible (generalized) $V_L$-module is equivalent to a $V_L$-module of the form $V_{L + \beta} \subset V_{L^{\circ}}$ for some $\beta \in L^{\circ}$ and any (generalized) $V_L$-module $W$ is equivalent to a direct sum of irreducible $V_L$-modules. i.e., \[W = \coprod_{\gamma_i
\in L^{\circ},\ i = 1, \dots, n} V_{\gamma_i + L},\] where
$\gamma_i$'s are arbitrary elements of $L^{\circ}$, and $n \in
\mathbb{N}$ (see \cite{D}, \cite{DLM}; cf. \cite{LL}).

Any strongly $M$-graded generalized $V_L$-module $W$ (in this example, all the generalized modules are modules) satisfies the assumption in Theorem \ref{main theorem} and the series (\ref{e3}), (\ref{e4}), (\ref{e5}) satisfies the expansions of the system of differential equations (\ref{eq7}) and (\ref{eq8}) in the regions $|z_1| > |z_2| > 0$, $|z_2| > |z_1 - z_2| > 0$, respectively.}
\end{exam}

\begin{exam}{\rm
Another family of strongly graded vertex algebra and modules comes from vertex algebras and modules associated to the abelian current Lie algebras. We recall the following construction from \cite{Y3}. Let $\mathfrak{h}$ be a finite dimensional abelian Lie algebra with a nondegenerate symmetric bilinear form $\langle\cdot, \cdot\rangle_{\mathfrak{h}}$. Let $\mathfrak{h}[t] = \mathfrak{h} \otimes \C[t]$. Then the current Lie algebra $\mathfrak{h}[t]$ is an abelian Lie algebra with an invariant symmetric bilinear form induced from $\langle\cdot, \cdot\rangle_{\mathfrak{h}}$. We form an affine Lie algebra
\[
\widehat{\mathfrak{h}[t]} = \mathfrak{h}[t]\otimes \C[s, s^{-1}] \oplus \C {\bf k}.
\]
It has graded subalgebras
\[
\widehat{\mathfrak{h}[t]}_{+} = \mathfrak{h}[t]\otimes s^{-1}\C[s^{-1}]
\]
and
\[
\widehat{\mathfrak{h}[t]}_{(\leq 0)} = \mathfrak{h}[t]\otimes \C[s] \oplus \C{\bf k}.
\]
Then the induced module
\[
M(l) = U(\widehat{\mathfrak{h}[t]})\otimes_{\widehat{\mathfrak{h}[t]}_{(\leq 0)}}\C{\bf 1} = S(\widehat{\mathfrak{h}[t]}_{+}) \otimes \C{\bf 1},
\]
where $\mathfrak{h}[t]\otimes \C[s]$ annihilates $\C{\bf 1}$ and ${\bf k}$ acts as a scalar multiplication by $l \in \C$, has a natural structure of quasi-conformal vertex algebra; see \cite{Y3}. For $i \in \N$, define
\[
M(l)^{(i)} = \mbox{span}\;\{h_1t^{i_1}(-n_1)\cdots h_kt^{i_k}(-n_k){\bf 1}\},
\]
where $h_j \in \mathfrak{h}$, $n_j \in \Z_{+}$ and $\sum_{j=1}^k i_j = i$ for $j = 1, \dots, k$. Thus $M(l)$ is equipped with a natural second grading by $\N$. Also for $n \in \N$, we have the doubly homogeneous subspace
\[
M(l)^{(i)}_{(n)} = \mbox{span}\;\{h_1t^{i_1}(-n_1)\cdots h_kt^{i_k}(-n_k){\bf 1}|\sum_{j=1}^k i_j = i, \sum_{j=1}^k n_j = n\},
\]
is finite dimensional, making $M(l)$ a strongly $\N$-graded quasi-conformal vertex algebra.

For $\lambda \in \mathfrak{h}^{*}$, let $\C_{\lambda}$ denote the one-dimensional $\mathfrak{h}$-module with $h \in \mathfrak{h}$ acts as $\lambda(h)$. For every $c \in \C$, we define the evaluation module $V(\lambda, c)$ as a vector space with the action given by
\[
(hf)\cdot v = f(c)\lambda(h) v\;\;\;h \in \mathfrak{h}\;f \in \C[t],\;v \in \C_{\lambda}.
\]
Then the induced module
\[
W(\lambda, 0, l) = U(\widehat{\mathfrak{h}[t]})\otimes_{\widehat{\mathfrak{h}[t]}_{(\leq 0)}}V(\lambda, 0) = S(\widehat{\mathfrak{h}[t]}_{+}) \otimes V(\lambda, 0)
\]
is a strongly $\N$-graded quasi-conformal module for $M(l)$.

For $\lambda \in \mathfrak{h}^{*}$, let $\Omega_{\lambda}$ denote a finite dimensional $\mathfrak{h}$-module such that
\[
(h -\lambda(h))^n v = 0
\]
for $h \in \mathfrak{h}$, $v \in \Omega_{\lambda}$ and some $n \in \Z_{+}$. For every $c \in \C$, we denote the evaluation $\widehat{\mathfrak{h}[t]}$-module induced from $\Omega_{\lambda}$ at $c$ by $\Omega(\lambda, c)$. We form the induced module
\[
G(\lambda, 0, l) = S(\widehat{\mathfrak{h}[t]}_{+}) \otimes \Omega(\lambda, 0).
\]
For simplicity, we assume $\mathfrak{h}$ is one-dimensional with nonzero vector $h$. If $h(0)^2$ admits a Jordan block with size greater than $1$ on $\Omega_{\lambda}$, then $G(\lambda, 0, l)$ will be a strongly $\N$-graded quasi-conformal generalized module for $M(l)$.

The strongly $\N$-graded module $W(\lambda, c, l)$ and generalized module $G(\lambda, 0, l)$ are $C_1$-cofinite with respect to $\N$, and the matrix elements of products and iterates of intertwining operators satisfy differential equations. See \cite{Y3} for the details.}
\end{exam}

\def\refname{\hfil{REFERENCES}}

\noindent {\small \sc Department of Mathematics, University of Notre Dame,
255 Hurley Building, Notre Dame, IN 46556-4618}
\vspace{1em}

\noindent {\em E-mail address}: jyang7@nd.edu


\begin{thebibliography}{fgtsth}
\bibitem[B1]{B1}
R. E. Borcherds, Monstrous moonshine and the monstrous Lie superalgebras, {\em Invent. Math.} {\bf 109} (1992), 405--444.

\bibitem[B2]{B2} R. E. Borcherds, Vertex algebras, Kac-Moody algebras, and the Monster, {\em Proc. Natl. Acad. Sci. USA} {\bf 83} (1986), 3068--3071.

\bibitem[D]{D} C. Dong, Vertex algebras associated with even
lattices, {\em J. of Algebra} {\bf 161} (1993), 245--265.


\bibitem[DLM]{DLM}
C.~Dong, H. Li and G. Mason, Regularity of rational vertex operator
algebras, {\it Adv. Math.} {\bf 132} (1997), 148--166.


\bibitem[FHL]{FHL}
I.~B. Frenkel, Y.-Z. Huang and J.~Lepowsky, {\em On axiomatic
approaches to vertex operator algebras and modules}, Mem. Amer.
Math. Soc. 104, Amer. Math. Soc., Providence, 1993, no. 494.

\bibitem[FLM]{FLM}
I.~B. Frenkel, J.~Lepowsky and A.~Meurman, {\em Vertex Operator
Algebras and the Monster}, Pure and Appl. Math., Vol. 134,  Academic
Press, Boston, 1988.

\bibitem[H]{H}
Y.-Z. Huang, Differential equations and intertwining operators, {\em Comm. Contemp.
Math}. {\bf 7} (2005), 375--400.

\bibitem[H1]{H1}
Y.-Z. Huang, Cofiniteness conditions, projective covers and the logarithmic tensor product theory, {\em J. Pure Appl. Alg.} {\bf 213} (2009), 458--475.

\bibitem[HLZ1]{HLZ1}
Y.-Z. Huang, J, Lepowsky and L. Zhang, Logarithmic tensor category theory for generalized modules for a conformal vertex algebra, I: Introduction and strongly graded algebras and their generalized modules, {\em Conformal Field Theories and Tensor Categories, Proceedings of a Workshop Held at Beijing International Center for Mathematics Research}, ed. C. Bai, J. Fuchs, Y.-Z. Huang, L. Kong, I. Runkel and C. Schweigert, Mathematical Lectures from Beijing University, Vol. 2, Springer, New York, 2014, 169--248.

\bibitem[HLZ7]{HLZ7}
Y.-Z. Huang, J, Lepowsky and L. Zhang, Logarithmic tensor category theory for generalized modules for a conformal vertex algebra, VII: Convergence and extension properties and applications to expansion for intertwining maps, arXiv: 1110.1929.

\bibitem[HLZ8]{HLZ8}
Y.-Z. Huang, J, Lepowsky and L. Zhang, Logarithmic tensor category theory for generalized modules for a conformal vertex algebra, VIII: Braided tensor category structure on categories of generalized modules for a conformal vertex algebra, arXiv: 1110.1931.

\bibitem[K]{K}
A. W. Knapp, {\em Representation Theory of Semisimple Groups,} Princeton University Press, Princeton, New Jersey, 1986.

\bibitem[LL]{LL}
J. Lepowsky and H. Li, {\em Introduction to Vertex Operator Algebras
and Their Representations}, Progress in Math., Vol. 227,
Birkh\"auser, Boston, 2003.

\bibitem[M1]{M1} A. Milas, Weak modules and logarithmic intertwining
operators for vertex operator algebras, in {\em Recent Developments
in Infinite-Dimensional Lie Algebras and Conformal Field Theory},
ed. S. Berman, P. Fendley, Y.-Z. Huang, K. Misra, and B. Parshall,
Contemp. Math., Vol. 297,  American Mathematical Society,
Providence, RI, 2002, 201--225.

\bibitem[M2]{M2} A. Milas, Logarithmic intertwining operators and vertex operators, {\em Comm. Math Phys.} {\bf 277} (2008), 497--529.

\bibitem[Y1]{Y1}
 J. Yang, Tensor products of strongly graded vertex algebras and their modules, {\em J. Pure Appl. Alg.} {\bf 217} (2013), 348--363.

\bibitem[Y2]{Y2}
J. Yang, Some results in the representation theory of strongly graded vertex algebras, Ph.D. thesis, Rutgers University, 2014.

\bibitem[Y3]{Y3}
J. Yang, Vertex algebra associated to the abelian current Lie algebras, arXiv: 1508.01923.
\end{thebibliography}
\end{document}